\documentclass[a4paper,11pt]{amsart}
\usepackage{times}
\usepackage[utf8]{inputenc}

\usepackage{amsmath,amssymb,amsthm,amsfonts,enumerate,url,fancyhdr,mathbbol,mathrsfs,color,hyperref}

\usepackage{pdfsync}

\usepackage{tikz,graphicx,pifont}
\usepackage{pgfplots}
\usepgfplotslibrary{polar}

\theoremstyle{plain}
\newtheorem{theorem}{Theorem}

\newtheorem{prop}[theorem]{Proposition}

\theoremstyle{definition}

\theoremstyle{remark}
\newtheorem{remark}[theorem]{Remark}

  \def\mG{\mathsf{G}}

  \def\mE{\mathsf{E}}

 \def\me{\mathsf{e}}

\def\XXint#1#2#3{{\setbox0=\hbox{$#1{#2#3}{\int}$ }
\vcenter{ \hbox{$#2#3$} }\kern-.45\wd0}}

\begin{document}

\title[The Cheeger constant of a quantum graph]{The Cheeger constant of a quantum graph}

\author{James B.\ Kennedy}
\address{James B.\ Kennedy, Institut f\"ur Analysis, Universit\"at Ulm, Helmholtzstr.\ 18, D-89081 Ulm, Germany and Institut f\"ur  Analysis, Dynamik und Modellierung, Universit\"at Stuttgart, Pfaffenwaldring 57, D-70569 Stuttgart, Germany}
\email{james.kennedy@mathematik.uni-stuttgart.de} 

\author{Delio Mugnolo}
\address{Delio Mugnolo, Lehrgebiet Analysis, Fakult\"at Mathematik und Informatik, Fern\-Universit\"at in Hagen, D-58084 Hagen, Germany}
\email{delio.mugnolo@fernuni-hagen.de} 

\begin{abstract}
We review the theory of Cheeger constants for graphs and quantum graphs and their present and envisaged applications.
\end{abstract}

\maketitle                   

The Laplacian matrix $\mathcal L$ of a graph $\mG$ (without loops)
has a long history, appearing and being rediscovered several times; be it as related to electrical circuits~\cite{Kir45}, to the discretisation of PDEs~\cite{Boo60}, to the theory of time-continuous Markov chains~\cite{Kat54} or to the formalism of Dirichlet forms~\cite{BeuDen59}.
In the 1970s Fiedler made the case for the study of the lowest non-zero eigenvalue $\lambda_1(\mG)$ of $\mathcal L$: since the multiplicity of $0$ as an eigenvalue of $\mathcal L$ is equal to the number of connected components of the underlying graph $\mG$, one may conjecture that if $0$ is a simple eigenvalue, then the smaller $\lambda_1(\mG)$, the closer the graph is to being disconnected. And indeed, the following important relation between $\lambda_1(\mG)$ and the edge connectivity $e(\mG)$ of $\mG$ (i.e., the minimal number of edges that have to be removed to make $\mG$ disconnected) was proved in~\cite[\S~4]{Fie73}. 

	\begin{prop}[Fiedler 1973]
Let $\mG$ be a connected graph on $V$ vertices. Then
\[
2e(\mG)\left(1-\cos\frac{\pi}{V}\right)\le \lambda_1(\mG)\le e(\mG)\ .
\]
\end{prop}
Thus, $e(\mG)$ and $\lambda_1(\mG)$ have the same asymptotic behaviour, although the scaling of $e(\mG)$ is sub-optimal as it penalises smaller graphs.
Adapting an idea developed in~\cite{Che70} for manifolds, several authors have studied since the beginning of the 1980s a renormalised version of $e(\mG)$: the \textit{Cheeger constant} $h(\mG)$ of $\mG$ is
\[
h(\mG):=\inf \frac{|\partial S|}{\min\{|S|, |S^C|\}}
\]
where $\inf$ is taken over all vertex sets $S$ and $\partial S$ is the set of all edges having exactly one endpoint in $S$ \cite{Dod84,AloMil85,Chu97}.

\begin{prop}[Dodziuk 1984, Alon--Milman 1985]\label{prop:alon}
Let $\mG$ be a connected graph of maximal degree $\deg_{\max}$. Then
\begin{equation}\label{eq:alon}
\frac{h^2(\mG)}{2\deg_{\max}(\mG)}\le \lambda_1(\mG)\le 2h(\mG)\ .
\end{equation}
\end{prop}
These estimates 
thus provide a variational relaxation of the NP-hard problem of determining $h(\mG)$~\cite{BuhHei09}.

Cheeger-type inequalities similar to~\eqref{eq:alon} hold for the Laplacian on \textit{quantum graphs} as well: recall that a quantum graph $\mathcal G$ is obtained from a graph $\mG$ by identifying each edge $\me$ with an interval $(0,\ell_\me)$. The \emph{standard Laplacian} on $\mathcal G$ is then a collection of second derivative operators on each edge, complemented with continuity and Kirchhoff (no flux loss) conditions in each edge \cite{BerKuc13,Mug14}. Nicaise introduced in~\cite{Nic87} a Cheeger-type constant for quantum graphs by
\[
h(\mathcal G):=\inf \frac{|\partial S|}{\min\{|S|, |S^C|\}}
\]
where $\inf$ is taken over all Lebesgue measurable open subsets $S$ of the quantum graph: $|\partial S|$ is the number of edges that depart from such $S$ and $|S|$ is its measure; \cite[Thm.~6.2]{DelRos16} characterises $h(\mathcal G)$ as the lowest non-zero eigenvalue of the $1$-Laplacian on $\mathcal G$.
 The straightforward estimate $h(\mathcal G)\ge \frac{2}{L}$ holds for all quantum graphs $\mathcal G$ of total length $L=\sum_{\me\in \mE}\ell_\me<\infty$; while the upper estimate $h (\mathcal G) \le \frac{2E}{L}$ -- with equality (among others) for flower graphs
with edges of equal length -- follows from~\cite[Thm.~6.2]{DelRos16} and (a straightforward extension of)~\cite[Lemma~2.3]{KenKurMal16}. Here $E$ is the number of \textit{essential edges}, i.e., the number of edges in $\mathcal G$ once vertices of degree 2 (irrelevant for the standard Laplacian) have been removed.

\begin{prop}\label{prop:nic}
Let $\mathcal G$ be a connected quantum graph with $E$ essential edges. Then the lowest non-zero eigenvalue  $\lambda_1(\mathcal G)$ of the standard Laplacian on $\mathcal G$ satisfies
\begin{equation}\label{eq:nic}
\max\left\{\frac{h^2(\mathcal G)}{4},\frac{\pi^2}{E^2}\frac{h^2(\mathcal G)}{4}\right\}\le \lambda_1(\mathcal G)\le \frac{\pi^2 E^2 h^2(\mathcal G)}{4}\ .\end{equation}
\end{prop}

The lower estimates in~\eqref{eq:nic} follow from~\cite[Théo.~3.2]{Nic87} and~\cite[Théo.~3.1]{Nic87} along with $h(\mathcal G)\le \frac{2E}{L}$; the upper estimate follows from~\cite[Thm.~4.2]{KenKurMal16} and $h(\mathcal G)\ge \frac{2}{L}$. We also mention the different but related upper estimate in~\cite[Thm.~1]{Kur13}. In analogy with a result obtained in~\cite{Par15} for convex subsets of $\mathbb R^2$, we conjecture that $\frac{\pi^2 h^2(\mathcal G)}{4}\le \lambda_1(\mathcal G)$.

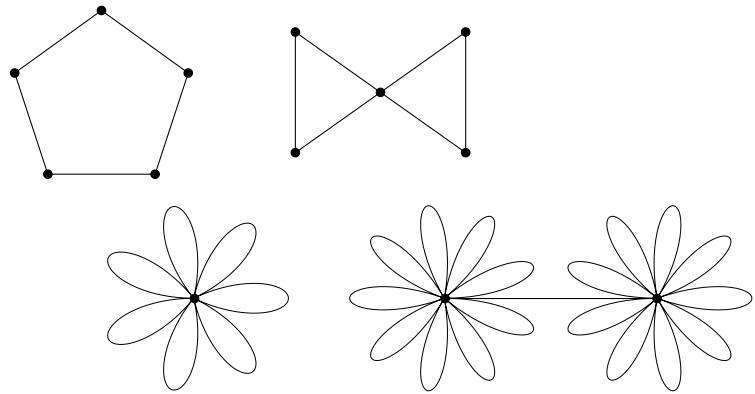
\begin{figure}
\begin{minipage}{3cm}
\begin{tikzpicture}[scale=0.8]
\draw (18:1.5cm) -- (90:1.5cm) -- (162:1.5cm) -- (234:1.5cm) -- (306:1.5cm) -- cycle;
\foreach \x in {18, 90, 162, 234, 306}{
\draw[fill] (\x:1.5cm) circle (2pt);}
\end{tikzpicture}
\end{minipage}
\begin{minipage}{.5cm}
$ $
\end{minipage}
\begin{minipage}{4.5cm}
\begin{tikzpicture}[scale=0.8]
\draw (4.6,-1) -- (4.6,1) -- (6,0) -- cycle;
\draw (6,0) -- (7.4,-1) -- (7.4,1) -- cycle;
\draw[fill] (4.6,-1) circle (2pt);
\draw[fill] (4.6,1) circle (2pt);
\draw[fill] (6,0) circle (2pt);
\draw[fill] (7.4,-1) circle (2pt);
\draw[fill] (7.4,1) circle (2pt);
\end{tikzpicture}
\end{minipage}\\[5pt]
\begin{minipage}{3.2cm}
\begin{tikzpicture}[scale=0.4]
  \begin{polaraxis}[grid=none, axis lines=none]
     \addplot[mark=none,domain=0:360,samples=300] { abs(cos(7*x/2))};
   \end{polaraxis}
\draw[fill] (3.45,3.42) circle (4pt);
 \end{tikzpicture}
\end{minipage}
\begin{minipage}{2.8cm}
\begin{tikzpicture}[scale=0.4]
\begin{scope}[xshift=4cm]
  \begin{polaraxis}[grid=none, axis lines=none]
     \addplot[mark=none,domain=0:360,samples=300]{abs(cos(9*x/2))/2};
   \end{polaraxis}
\end{scope}
\draw (0.45,3.42) -- (7.42,3.42);
\draw[fill] (0.45,3.42) circle (4pt);
\draw[fill] (7.42,3.42) circle (4pt);
\begin{scope}[xshift=-3cm]
  \begin{polaraxis}[grid=none, axis lines=none]
     \addplot[mark=none,domain=0:360,samples=300] { abs(sin(9*x/2))};
   \end{polaraxis}
\end{scope}
 \end{tikzpicture}
\end{minipage}
\caption{Four quantum graphs: a cycle, a butterfly, a flower and a symmetric flower dumbbell}\label{fig:1}
\end{figure}

\begin{remark}
1) If $\mathcal G$ is an interval, then $\lambda_1(\mathcal G)=\frac{\pi^2}{L^2}=\frac{\pi^2 h^2(\mathcal G)}{4}$. If $\mathcal G$ is 
a flower, then $\lambda_1(\mathcal G)=\frac{\pi^2E^2}{L^2}=\frac{\pi^2 h^2(\mathcal G)}{4}$. Unfortunately, the dependence on $E$ cannot in general be dropped in the upper estimate in~\eqref{eq:nic}: \textit{symmetric flower dumbbells} (see Figure~\ref{fig:1}) obviously have Cheeger constant $\frac{2}{L}$, as the optimal Cheeger set $S$ is obtained by just cutting $\mathcal G$ in the middle. 
At the same time, by adding more and more petals and simultaneously shortening all of them while making the handle shorter and shorter, one can produce symmetric flower dumbbells with same total length but arbitrarily high $\lambda_1(\mathcal G)$.

2) While we do not know whether the upper estimate in~\eqref{eq:nic} is sharp, symmetric flower dumbbells with $E=2m+1$ edges satisfy
\begin{equation}\label{eq:nic-3}
\lambda_1(\mathcal G)\approx {\pi^2m^2}= \frac{\pi^2 (E-1)^2 h^2(\mathcal G)}{4} \ ,
\end{equation}
which is the corresponding value of $\lambda_1$ for a flower with $E-1$ edges, provided the symmetric flower dumbbell's handle is arbitrarily short.
\end{remark}

The main fascinating feature of the Cheeger constant of quantum graphs is its hybrid nature, partly combinatorial and partly metric (its numerator and denominator, respectively), in sharp contrast to its counterparts for manifolds and graphs. But
is it meaningful at all to consider the Cheeger constant of a quantum graph?
From the point of view of theoretical computer science the lowest non-zero eigenvalue is an elementary object that can be easily determined by variational methods and can in turn help to estimate the Cheeger constant -- the \emph{really} interesting quantity, for the purpose of machine learning. 

We maintain that quantum graphs are not  unnecessarily complicated gadgets, but rather useful tools delivering additional information. 
As an example, let us consider the first two quantum graphs in Figure~\ref{fig:1}, each of whose intervals is assumed to have unit length.
One sees that the Cheeger constant of the cycle is $\frac{4}{5}$, while the  butterfly has Cheeger constant $\frac{2}{3}$. On the other hand, both underlying \textit{discrete} graphs have Cheeger constant 1. We argue that the information yielded by $h(\mG)$ may in critical cases be \emph{complemented} by $h(\mathcal G)$, upon turning a graph $\mG$ into a quantum graph $\mathcal G$ with edges of unit length, whenever the \emph{interaction-based} description offered by a quantum graph is as relevant as the \emph{agent-based} description offered by a graph.

\begin{thebibliography}{KKMM16}

\bibitem[AM85]{AloMil85}
N.~Alon and V.~D. Milman.
\newblock {$\lambda_1,$} isoperimetric inequalities for graphs, and
  superconcentrators.
\newblock {\em J.\ Combin.\ Theory Ser.\ B}, 38:73--88, 1985.

\bibitem[BD59]{BeuDen59}
A.\ Beurling and J.~Deny.
\newblock {Dirichlet spaces}.
\newblock {\em Proc.\ Natl.\ Acad.\ Sci.\ USA}, 45:208--215, 1959.

\bibitem[BH09]{BuhHei09}
T.~Bühler and M.~Hein.
\newblock Spectral clustering based on the graph $p$-{L}aplacian.
\newblock In {\em Proc.\ 26th Annual Int.\ Conf.\ Mach.\ Learning}, pages
  81--88, New York, 2009. ACM.

\bibitem[BK13]{BerKuc13}
G.\ Berkolaiko and P.~Kuchment.
\newblock {\em {Introduction to Quantum Graphs}}, volume 186 of {\em Math.\
  Surveys and Monographs}.
\newblock Amer.\ Math.\ Soc., Providence, RI, 2013.

\bibitem[Boo60]{Boo60}
G.~Boole.
\newblock {\em {A Treatise on the Calculus of Finite Differences}}.
\newblock Macmillan, Cambridge, 1860.

\bibitem[Che70]{Che70}
J.\ Cheeger.
\newblock A lower bound for the smallest eigenvalue of the {L}aplacian.
\newblock In R.C. Gunning, editor, {\em Problems in Analysis}, pages 195--199,
  Princeton, NJ, 1970. Princeton Univ.\ Press.

\bibitem[Chu97]{Chu97}
F.R.K. Chung.
\newblock {\em Spectral Graph Theory}, volume~92 of {\em Reg.\ Conf.\ Series
  Math.}
\newblock Amer.\ Math.\ Soc., Providence, RI, 1997.

\bibitem[Dod84]{Dod84}
J.~Dodziuk.
\newblock Difference equations, isoperimetric inequality and transience of
  certain random walks.
\newblock {\em Trans.\ Amer.\ Math.\ Soc.}, 284:787--794, 1984.

\bibitem[DPR16]{DelRos16}
L.M. Del~Pezzo and J.D. Rossi.
\newblock The first eigenvalue of the $p$-{L}aplacian on quantum graphs.
\newblock {\em Analysis and Math.\ Phys.}, DOI:10.1007/s13324-016-0123-y, 2016.

\bibitem[Fie73]{Fie73}
M.~Fiedler.
\newblock Algebraic connectivity of graphs.
\newblock {\em Czech.\ Math.\ J.}, 23:298--305, 1973.

\bibitem[Kat54]{Kat54}
T.~Kato.
\newblock {On the semi-groups generated by Kolmogoroff's differential
  equations}.
\newblock {\em J.\ Math.\ Soc.\ Jap.}, 6:1--15, 1954.

\bibitem[Kir45]{Kir45}
G.~Kirchhoff.
\newblock {Ueber den Durchgang eines elektrischen Stromes durch eine Ebene,
  insbesondere durch eine kreisf\"ormige}.
\newblock {\em Ann.\ Physik}, 140:497--514, 1845.

\bibitem[KKMM16]{KenKurMal16}
J.B. Kennedy, P.~Kurasov, G.~Malenová, and D.~Mugnolo.
\newblock On the spectral gap of a quantum graph.
\newblock {\em Ann.\ Henri Poincar\'e A}, DOI:10.1007/s00023-016-0460-2, 2016.

\bibitem[Kur13]{Kur13}
P.~Kurasov.
\newblock On the spectral gap for {L}aplacians on metric graphs.
\newblock {\em Acta Phys.\ Pol.\ A}, 124:1060--1062, 2013.

\bibitem[Mug14]{Mug14}
D.~Mugnolo.
\newblock {\em {Semigroup Methods for Evolution Equations on Networks}}.
\newblock Springer-Verlag, Berlin, 2014.

\bibitem[Nic87]{Nic87}
S.\ Nicaise.
\newblock Spectre des r{\'e}seaux topologiques finis.
\newblock {\em Bull.\ Sci.\ Math., II.\ S{\'e}r.}, 111:401--413, 1987.

\bibitem[Par15]{Par15}
E.~Parini.
\newblock Reverse cheeger inequality for planar convex sets.
\newblock {\em arXiv:1501.04520}, 2015.

\end{thebibliography}
\end{document}